\documentclass[12pt]{article}
\usepackage{geometry}

\usepackage{amsfonts}
\usepackage{amsmath}
\usepackage{amssymb}
\usepackage{amsthm}

\newcommand{\ceil}[1]{{\lceil #1 \rceil}}
\newcommand{\floor}[1]{{\lfloor #1 \rfloor}}

\newcommand{\dist}{{\rm dist}}

\newcommand{\Set}[1]{\left\{ #1 \right\}}
\newcommand{\Z}{{\mathbb Z}}

\newcommand{\cx}{{\mathbb C}}

\newcommand{\fracpart}[1]{\langle #1 \rangle}

\newcommand{\half}{\frac{1}{2}}

\newcommand{\rat}{{\mathbb Q}}
\newcommand{\real}{{\mathbb R}}

\newcommand{\rme}{{\rm e}}
\newcommand{\rmi}{{\sqrt{-1}}}
\newcommand{\set}[1]{\{ #1 \}}

\newcommand{\zahl}{{\mathbb Z}}

\newtheorem{theorem}{Theorem}

\newtheorem{definition}[theorem]{Definition}
\newtheorem{lemma}[theorem]{Lemma}
\newtheorem{proposition}[theorem]{Proposition}

\title{Spiral Delone sets in relative metric}
\author{Yoshikazu Yamagishi\thanks{Ryukoku University. yg@rins.ryukoku.ac.jp}}

\begin{document}

\maketitle

\begin{abstract}
A general Archimedean spiral lattice is a Delone set in the relative
distance if and only if its rotation angle is badly
approximable.
\end{abstract}

\medskip


\section{Introduction}

Phyllotaxis is a subject in biology which is related to
crystallography and mathematics.  It studies the arrangements of
leaves, seeds, and other organs of plants.  In the 19th century,
botanical spiral patterns were considered as living crystals by
A.~Bravais and L.~Bravais, who assumed that their rotation angles are
irrational numbers.

The arrangement of seeds in the sunflower head
is modeled by the Bernoulli spiral lattice $\set{r^n \rme^{2\pi n \tau\rmi}
  \mid n \in \zahl}$, $\tau = (1+\sqrt{5})/2$, $0 < r < 1$, and the
general Archimedean spiral lattice $\set{n^\alpha\rme^{2\pi n\tau\rmi}
  \mid n \in \zahl_{\ge0}}$, $\alpha>0$ (the case $\alpha=1/2$ is
called Fermat, $\alpha=1$ Archimedean).

Let 
\begin{equation}
\label{40}
\Gamma_{\alpha,\theta} = \set{ n^\alpha \rme^{2\pi n\theta\rmi} \mid n \in \zahl_{\ge0}},
\quad \alpha >0, \ \theta \in \real.
\end{equation}
Akiyama \cite{akiyama} showed that (i) if $\alpha > \half$ then
$\Gamma_{\alpha,\theta}$ is not relatively dense, (ii) if $\alpha <
\half$ then $\Gamma_{\alpha,\theta}$ is not uniformly discrete, (iii)
if $\alpha = \half$ and $\theta$ is irrational, then
$\Gamma_{\half,\theta}$ is relatively dense $\Leftrightarrow$
$\Gamma_{\half,\theta}$ is uniformly discrete $\Leftrightarrow$
$\theta$ is badly approximable.  His proof is based on the Three Gap
Theorem.  A Delone set is a set which is both relatively dense and
uniformly discrete.  So a Fermat spiral lattice is a Delone set if and
only if $\theta$ is badly approximable.

This paper addresses the relative denseness and uniform discreteness
of $\Gamma_{\alpha,\theta}$ with respect to the relative distance
\[ d_\beta(z,w):= \frac{|z-w|}{|z|^\beta + |w|^\beta}, 
 \quad -\infty < \beta < 1.
\]
It is known (\cite{hasto2002}) that if $\half \le \beta \le 1$, the
function $d_\beta$ satisfies the triangle inequality, so that it is a
metric on $\cx$.  In this paper we define the notions of asymptotic
$\beta$-relative denseness and asymptotic $\beta$-uniform discreteness
with respect to the relative distance $d_\beta$.  We show that

\begin{lemma}
\label{41}
If $\Gamma_{\alpha,\theta}$ is asymptotically $\beta$-relatively
dense, then $1 + 2\alpha\beta \ge 2\alpha$.
\end{lemma}

\begin{lemma}
\label{42}
If $\Gamma_{\alpha,\theta}$ is asymptotically $\beta$-uniformly
discrete, then $1 + 2\alpha\beta \le 2\alpha$.
\end{lemma}

\begin{lemma}
\label{43}
Suppose that $1 + 2\alpha\beta = 2\alpha$.  If $\Gamma_{\alpha,\theta}$ is
asymptotically $\beta$-relatively dense or asymptotically
$\beta$-uniformly discrete, then $\theta \not\in\rat$.
\end{lemma}

\begin{theorem}
\label{44}
Suppose that $1 + 2\alpha\beta = 2\alpha$.
Then the following statements are mutually equivalent.
\begin{itemize}
\item[(i)] $\Gamma_{\alpha,\theta}$ is asymptotically
  $\beta$-relatively dense.
\item[(ii)] $\Gamma_{\alpha,\theta}$ is asymptotically
  $\beta$-uniformly discrete.
\item[(iii)] $\theta$ is badly approximable.
\end{itemize}
\end{theorem}

A main idea of the proof of Theorem \ref{44} is the approximation by
(linear) lattices.  Denote by $- \half \le \fracpart{x} < \half$ be
the fractional part of $x \in \real$, where $x - \fracpart{x} \in
\zahl$.  By Taylor's theorem, there exists a constant $C>0$ depending
only on $\alpha$, such that
\begin{align*}
(\nu + k)^\alpha \rme^{2\pi k\theta \rmi} 
&= \nu^\alpha \left( 1 + \frac{k}{\nu} \right)^\alpha \rme^{2\pi k\theta\rmi} 
\notag \\
&= \nu^\alpha \left( 1 + \frac{\alpha k}{\nu} + 2\pi \fracpart{\theta k} \rmi
 \right) + \epsilon
\notag \\
&= \nu^\alpha + \nu^{\alpha-\half}
 \left( \frac{\alpha k}{\sqrt{\nu}} + 2\pi \fracpart{\theta k} \sqrt{\nu} \rmi
  \right) + \epsilon, 
\end{align*}
where
\[
|\epsilon| \le C \nu^\alpha  
\left( \left| \frac{\alpha k}{\nu} \right|^2 + |2\pi\fracpart{\theta k}|^2
\right).
\]

Before proving Theorem \ref{44}, we consider the family of linear lattices
\begin{equation}
\label{4}
\Lambda(t) 
= \Set{ \left( (m\theta - n)\sqrt{t}, \frac{m}{\sqrt{t}} \right)
  \mid m,n\in\zahl}, 
 \quad t \ge 1.
\end{equation}
It can be regarded as a scenery flow (\cite{arnoux-fisher}) of the
lattice $\Lambda = \set{ (m\theta-n, m) \mid m,n\in\zahl}$.  We define
the notions of relative denseness and uniform discreteness of the
family $\set{\Lambda(t)}_{t\ge1}$, and show the following Proposition.

\begin{proposition}
\label{7}
Let $\theta \in \real\setminus\rat$.  Let $\set{\Lambda(t)}_{t\ge1}$
be a family of lattices defined in (\ref{4}).  Then the following
conditions are mutually equivalent.
\begin{itemize}
\item[(i)] $\set{\Lambda(t)}_{t\ge1}$ is relatively dense.
\item[(ii)] $\set{\Lambda(t)}_{t\ge1}$ is uniformly discrete.
\item[(iii)] $\theta$ is badly approximable.
\end{itemize}
\end{proposition}

Proposition \ref{7} works as a linear prototype of Theorem \ref{44}.
Note that a key tool in the proof of Proposition \ref{7} is Richards'
formula, which is known in the phyllotaxis theory \cite{jean}.

Section 2 defines the notions of asymptotically $\beta$-relative
denseness and asymptotically $\beta$-uniform discreteness, and show
Lemmas \ref{41}-\ref{43}.  Setion 3 prepares notations in the
continued fraction expansions and rational approximations.  Section 4
proves Proposition \ref{7}, by using Richards' formula.  Section 5
proves Theorem \ref{44}.

See \cite{jean,adler} for the history of the study of phyllotaxis, and
\cite{pennybacker,barabe} for recent surveys.  In \cite{arxiv2020}, it
was shown that the area of Voronoi cells for a general Archimedean
spiral lattice has a convergence under some scale normalization.  In
\cite{jpa2018}, the combinatorial structures of the grains and grain
boundaries of the Voronoi tessellations for Archimedean spiral
lattices were described.  In \cite{physd2017}, it was shown that, in
the family of Bernoulli spiral lattices, the bifurcation diagram of
Voronoi tessellations is a dual graph of the bifurcation diagram of
circle packings, by using the relative metric
$d(z,w)=|z-w|/(|z|+|w|)$.  Marklof \cite{marklof} showed that the
point set $\set{ \sqrt{n} \rme^{2\pi\theta\sqrt{n}\rmi} \mid n \in
  \zahl_{\ge0}}$ is a Delone set for any $\theta > 0$.

\section{Delone sets in relative metric}
\label{27}

In this section we define the notions of asymptotically
$\beta$-relative denseness and asymptotically $\beta$-uniform
discreteness, and show Lemmas \ref{41}-\ref{43}.

\begin{lemma}
\label{15}
Let $r', r > 0$, $0 < \beta < 1$, $z, \zeta \in \cx$.  Suppose that $r
< r' < 2r$.  If $|z - \zeta| < r |z|^\beta$ and $|z| \ge M := \left(
\frac{rr'}{2(r'-r)}\right)^{1/(1-\beta)}$, then we have $|z - \zeta| <
\frac{r'}{2}|z|^\beta + \frac{r'}{2}|\zeta|^\beta$.
\end{lemma}
\begin{proof}
If $|\zeta| \ge |z|$, we have $|z - \zeta| < r|z|^\beta < r' |z|^\beta
\le \frac{r'}{2}|z|^\beta + \frac{r'}{2}|\zeta|^\beta$.  So we assume
that $|z| > |\zeta|$.  Since $|z| \ge M =
(\frac{rr'}{2(r'-r)})^{1/(1-\beta)}$, we have $|z| \ge
\frac{rr'}{2(r'-r)} |z|^\beta$, and $r|z|^\beta \le
\frac{2(r'-r)}{r'}|z|$.  Hence $|z| - |\zeta| \le |z - \zeta| < r
|z|^\beta \le \frac{2(r'-r)}{r'}|z|$.  This implies that $0 <
\frac{2r-r'}{r'} < \frac{|\zeta|}{|z|} < 1$, and $\frac{2r-r'}{r'} <
\frac{|\zeta|}{|z|} < \frac{|\zeta|^\beta}{|z|^\beta}$ since $0<\beta
< 1$.  So we obtain $(2r-r')|z|^\beta < r'|\zeta|^\beta$, and $|z -
\zeta| < r |z|^\beta < \frac{r'}{2}(|z|^\beta + |\zeta|^\beta)$.
\end{proof}

\begin{lemma}
\label{2}
Let $r' > r > 0$, $0 < \beta < 1$, $z, \zeta \in \cx$.
Suppose that
$|z| \ge M := \left( \frac{rr'}{r'-r}\right)^{1/(1-\beta)}$.
If either  
\begin{equation}
\label{17}
|z - \zeta| < \frac{r}{2}|z|^\beta + \frac{r}{2}|\zeta|^\beta
\end{equation}
or
\begin{equation}
\label{16}
 |z - \zeta| < r |z|^\beta,
\end{equation}
then we have $|z - \zeta| < r'|\zeta|^\beta$.
\end{lemma}
\begin{proof}
If $|\zeta| \ge |z|$, we have
$r|z|^\beta \le \frac{r}{2}|z|^\beta + \frac{r}{2}|\zeta|^\beta
\le r|\zeta|^\beta < r' |\zeta|^\beta$.
So we assume that $|z| > |\zeta|$.
Then we have
$\frac{r}{2}|z|^\beta + \frac{r}{2}|\zeta|^\beta
  < r |z|^\beta$,
so the assumption (\ref{17}) implies (\ref{16}).
Since $|z| \ge M = (\frac{rr'}{r'-r})^{1/(1-\beta)}$,
we have $|z| \ge \frac{rr'}{r'-r} |z|^\beta$,
and $r|z|^\beta \le \frac{r'-r}{r'}|z|$.
Hence
$|z| - |\zeta|
\le |z - \zeta|
 <  r |z|^\beta \le \frac{r'-r}{r'}|z|$.
This implies that
$\frac{r}{r'} < \frac{|\zeta|}{|z|} < 1$,
and $\frac{r}{r'} < \frac{|\zeta|}{|z|} < \frac{|\zeta|^\beta}{|z|^\beta}$
since $0<\beta < 1$.
So we obtain $|z - \zeta| < r |z|^\beta < r'|\zeta|^\beta$.
\end{proof}

\begin{lemma}
\label{18}
If $\alpha >0$ and $0 < x \le \frac{1}{\alpha^2}$, then
$(1+x)^\alpha < 1 + 2\alpha x$.
\end{lemma}
\begin{proof}
If $0 < \alpha \le 1$, then
$(1+x)^\alpha \le 1+\alpha x < 1 + 2\alpha x$.  
So suppose that $\alpha >1$.
Since the function $t \mapsto \rme^t$ is convex, we have
$\frac{\rme^{1/\alpha} - 1}{1/\alpha} \le \rme^1 - 1 \le 2$.
Since the function $x \mapsto (1+x)^\alpha$ is convex, we have
$\frac{(1+x)^\alpha - 1}{x}
\le \frac{(1+\frac{1}{\alpha^2})^\alpha - 1}{1/\alpha^2}
\le \frac{\rme^{1/\alpha} - 1}{1/\alpha^2}
\le 2\alpha$, which completes the proof.
\end{proof}

\begin{lemma}
If $\delta>0$ and $0 < h < \half$, then $(1-h)^\delta > 1 - 2\delta h$.
\end{lemma}
\begin{proof}
If $\delta \ge 1$, then $(1-h)^\delta \ge 1 - \delta h \ge 1 - 2\delta
h$.  Suppose that $0 < \delta < 1$.  By the Mean Value Theorem, we
have $(1 - (1-h)^\delta)/h = \delta (1-\epsilon h)^{\delta - 1}$ for
some $0 < \epsilon < 1$, and $(1 - \epsilon h)^{\delta - 1} < (1 -
h)^{-1} \le 2$ for $0 < h < \half$.
\end{proof}

\begin{lemma}
\label{19}
Let $r' > r > 0$, $-\infty < \beta < 0$, $\delta = -\beta > 0$,
$z,\zeta \in \cx$.  Suppose that $|z| \ge M: = \max\set{1, r\delta^2,
  \frac{2r r'\delta}{r'-r}, 2r}$ and $|\zeta|\ge1$.  Then the
following statements hold.
\begin{enumerate}
\item[(1)]
If $|z - \zeta| < r|z|^\beta$,
then we have $|z - \zeta| < r' |\zeta|^\beta$
and $|z - \zeta| < \frac{r'}{2}(|z|^\beta + |\zeta|^\beta)$.
\item[(2)]
If $|z - \zeta| < \frac{r}{2}(|z|^\beta + |\zeta|^\beta)$,
then $|z - \zeta| < r' |\zeta|^\beta$.
\item[(3)]
If $|z - \zeta| < r |\zeta|^\beta$, 
then $|z - \zeta| < r' |z|^\beta$.
\end{enumerate}
\end{lemma}
\begin{proof}
(1). \ 
If $|\zeta| \le |z|$, we have
$|z - \zeta| < r |z|^{-\delta} < \frac{r'}{2}(|z|^{-\delta} + |\zeta|^{-\delta})
 < r' |\zeta|^{-\delta}$. 
So assume that $|\zeta| > |z|$.
Since $|\zeta| - |z| \le |z - \zeta| < r |z|^{-\delta} \le r$,
we have $|\frac{\zeta}{z}| < 1 + \frac{r}{|z|}
< 1 + \frac{r}{M}$.
Since $\frac{r}{M} \le \frac{1}{\delta^2}$,
Lemma~\ref{18} applies to see that
\[
\frac{|\zeta|^\delta}{|z|^\delta}
< \left( 1 + \frac{r}{M} \right)^\delta
\le 1 + \frac{2r\delta}{M}
\le \frac{r'}{r},
\]
so we obtain $|z - \zeta| < r|z|^{-\delta} 
 < r' |\zeta|^{-\delta} < 
\frac{r'}{2}(|z|^{-\delta} + |\zeta|^{-\delta})$.

(2).\
If $|\zeta| \le |z|$, we have
$|z - \zeta| < \frac{r}{2}(|z|^{-\delta} + |\zeta|^{-\delta})
< r' |\zeta|^{-\delta}$. 
If $|\zeta| > |z|$, we have
$|z - \zeta|
< \frac{r}{2}(|z|^{-\delta} + |\zeta|^{-\delta}) < r|z|^{-\delta}$,
so (1) applies to see that  $|z - \zeta| < r' |\zeta|^{-\delta}$.

(3).\ 
If $|\zeta| \ge |z|$, we have
$|z - \zeta| < r |\zeta|^{-\delta} < r' |z|^{-\delta}$. 
Next assume that $1 \le |\zeta| < |z|$. 
Since $|z| - |\zeta| \le |z - \zeta| < r |\zeta|^{-\delta} \le r$,
we have $| \frac{\zeta}{z}| \ge 1 - \frac{r}{|z|} \ge 1 - \frac{r}{M}$.
So
\[ \left| \frac{\zeta}{z} \right|^\delta 
 \ge \left( 1 - \frac{r}{M} \right)^\delta
 \ge 1 - \frac{2r\delta}{M}
 \ge \frac{r}{r'},
\]
and we obtain $|z - \zeta| < r |\zeta|^{-\delta} \le r' |z|^{-\delta}$.
\end{proof}

\begin{lemma}
\label{3}
Let $- \infty < \beta < 1$, $r>0$.  Let $\Gamma \subset \cx$.
Let $\Gamma' = \set{\gamma \in \Gamma \mid |\gamma|\ge1}$.
The following conditions are mutually equivalent.
\begin{itemize}
\item[(i)] For any $r_1>r$ there exists $M_1 >0$ such that for any $z
  \in \cx$ with $|z| \ge M_1$, there exists $\zeta \in \Gamma'$ such
  that $|z-\zeta| < r_1 |z|^\beta$.
\item[(ii)] For any $r_2 >r$ there exists $M_2 >0$ such that for any
  $z \in \cx$ with $|z| \ge M_2$, there exists $\zeta \in \Gamma'$
  such that $|z-\zeta| < r_2 |\zeta|^\beta$.
\item[(iii)] For any $r_3 > r$ there exists $M_3 >0$ such that for any
  $z \in \cx$ with $|z| \ge M_3$, there exists $\zeta \in \Gamma'$
  such that $|z-\zeta| < \frac{r_3}{2} (|z|^\beta + |\zeta|^\beta)$.
\end{itemize}
\end{lemma}

\begin{proof}
(i)$\Rightarrow$(ii).  The case $\beta = 0$ is trivial, so suppose
  that $-\infty < \beta < 0$ or $0 < \beta < 1$.  Let $r_2 >r$.  Take
  $r_1$ such that $r < r_1 < r_2$.  There exists $M_1>0$ such that for
  each $z \in \cx$ with $|z|\ge M_1$, there exists $\zeta \in \Gamma$
  such that $|z - \zeta| < r_1 |z|^\beta$.  Let $M_2 =\max\Set{M_1,
    \left( \frac{r_1 r_2}{r_2-r_1}\right)^{1/(1-\beta)}}$ if $0 <
  \beta < 1$, or let $M_2 = \max\Set{M_1, r_1 \delta^2, \frac{2r_1
      r_2\delta}{r_2 - r_1}}$ if $-\infty < \beta < 0$, where $\delta
  = - \beta > 0$.  If $|z| \ge M_2$, then we have $|z - \zeta| < r_2
  |\zeta|^\beta$ by Lemmas~\ref{2}, \ref{19}.

(ii)$\Leftrightarrow$(iii)$\Rightarrow$(i).
All the other arguments are given in a similar way.
\end{proof}

A point set $\Gamma \subset \cx$ is called locally finite if
$\#\set{\zeta \in \Gamma \mid |\zeta|<r} < \infty$ for any $r>0$.

\begin{lemma}
\label{26}
Let $- \infty < \beta < 1$, $s>0$.  Suppose that a point set $\Gamma
\subset \cx$ is locally finite.  Let $\Gamma' = \set{\gamma\in \Gamma
  \mid |\gamma|\ge 1}$.  The following conditions are mutually
equivalent.
\begin{itemize}
\item[(i)] For any $0 < s_1 < s$, there exists $M_1>0$ such that for
  any $z \in \cx$ with $|z| \ge M_1$, we have $\#\set{\zeta \in
    \Gamma' \mid |z-\zeta| < s_1 |z|^\beta} \le 1$.
\item[(ii)] For any $0 < s_2 < s$, there exists $M_2>0$ such that for
  any $z \in \cx$ with $|z| \ge M_2$, we have $\#\set{\zeta \in
    \Gamma' \mid |z-\zeta| < s_2 |\zeta|^\beta} \le 1$.
\item[(iii)] For any $0 < s_3 < s$, there exists $M_3>0$ such that for
  any $z \in \cx$ with $|z| \ge M_3$, we have $\#\set{\zeta \in
    \Gamma' \mid |z-\zeta| < \frac{s_3}{2} (|z|^\beta +
    |\zeta|^\beta)} \le 1$.
\end{itemize}
\end{lemma}

\begin{proof}
(i)$\Rightarrow$(ii).  The case $\beta=0$ is trivial, so assume that
  $-\infty < \beta < 0$ or $0 < \beta < 1$.  Suppose that (ii) does
  not hold.  There exists $s_2 < s$ and a sequence $z_i$, $i \in
  \zahl_{>0}$, such that $\lim_{i\to\infty}|z_i|=+\infty$ and
  $\#\set{\zeta \in \Gamma' \mid |z_i - \zeta| <
    s_2|\zeta|^\beta}\ge2$ for each $i \in \zahl_{>0}$.  Take $s_1$
  such that $s_2 < s_1 < s$.  Let $M=\left( \frac{s_1 s_2}{s_1 -
    s_2}\right)^{1/(1-\beta)}$ if $0 < \beta < 1$, or let $M =
  \max\Set{s_2\delta^2, \frac{2s_1 s_2\delta}{s_1 - s_2}, 2s_2}$ if
  $-\infty < \beta < 0$, where $\delta = - \beta > 0$.  If $|z_i|\ge
  M$, then we obtain $\#\set{\zeta \in \Gamma' \mid |z_i - \zeta| <
    s_1|z|^\beta}\ge2$ by Lemmas~\ref{2}, \ref{19}.  So (i) does not
  hold.

(ii)$\Leftrightarrow$(iii)$\Rightarrow$(i).
All the other arguments are given in a similar way.
\end{proof}

\begin{definition}
We say that a point set $\Gamma \subset \cx$ is asymptotically
$\beta$-relatively dense if there exists $r>0$ that satisfies one (and
hence all) of the conditions in Lemma \ref{3}.
\end{definition}

\begin{definition}
We say that $\Gamma \subset \cx$ is asymptotically $\beta$-uniformly
discrete if there exists $r>0$ that satisfies one (and hence all) of
the conditions in Lemma \ref{26}.  $\Gamma$ is called an
asymptotically $\beta$-Delone set if it is both asymptotically
$\beta$-relatively dense and asymptotically $\beta$-uniformly
discrete.
\end{definition}

Now we prove Lemmas \ref{41}-\ref{43}.
Let $B(z,r) = \set{\zeta \in \cx \mid |\zeta - z| < r}$,
$z \in \cx$, $r >0$, be an open disk.

\begin{proof}[Proof of Lemma \ref{41}]
If $\Gamma_{\alpha,\theta}$ is asymptotically $\beta$-relatively
dense, then there exist $r>0$ and $m_0 \in \zahl_{>0}$ such that the
region $\set{\zeta \in \cx \mid \zeta \ge m_0^\alpha + r
  m_0^{\alpha\beta}}$ is covered by the family of disks
$\set{B(n^\alpha \rme^{2\pi n\theta\rmi}, r n^{\alpha\beta}) \mid n
  \in \zahl, n \ge m_0}$.

First suppose that $0 \le \beta < 1$.
Then we have
\[ \sum_{m_0 < n \le m} \pi (r n^{\alpha\beta})^2
\ge \pi (m^\alpha - r m^{\alpha\beta})^2 - \pi (m_0^\alpha + r m_0^{\alpha\beta})^2
\]
for any $m > m_0$, so
\[ (m - m_0) \pi (r m^{\alpha\beta})^2 \ge 
\pi (m^\alpha - r m^{\alpha\beta})^2 - \pi (m_0^\alpha + r m_0^{\alpha\beta})^2
\]
By taking $m \to \infty$, we obtain $1 + 2\alpha \beta \ge 2 \alpha$.

If $-\infty < \beta < 0$, we have
\[ \sum_{m < n \le 2m} \pi (r n^{\alpha\beta})^2
\ge \pi ((2m)^\alpha - r (2m)^{\alpha\beta})^2 - \pi (m^\alpha + r m^{\alpha\beta})^2
\]
for any $m \ge m_0$.  We may assume $m$ is so large that $(2m)^\alpha
- r (2m)^{\alpha\beta} \ge (\frac{3}{2})^\alpha (m^\alpha + r
m^{\alpha\beta})$.  Then
\[ m\pi (r m^{\alpha\beta})^2 
 \ge \pi \left( \left(\frac{3}{2} \right)^\alpha - 1 \right)
      (m^\alpha + r m^{\alpha\beta})^2.
\]
By taking $m \to \infty$, we obtain $1 + 2\alpha \beta \ge 2\alpha$.
\end{proof}

\begin{proof}[Proof of Lemma \ref{42}]
If $\Gamma_{\alpha,\theta}$ is asymptotically $\beta$-uniformly
discrete, then there exist $r>0$ and $m_0 \in \zahl_{>0}$ such that
the family of disks $\set{B(n^\alpha \rme^{2\pi n\theta\rmi}, r
  n^{\alpha\beta}) \mid n\in\zahl, n > m_0}$, is distinct.  We assume
that $m_0$ is sufficiently large and $m_0^\alpha - r m_0^{\alpha\beta}
> 0$.

If $0 \le \beta < 1$,
we have
\[ \sum_{m < n \le 2m} \pi (r n^{\alpha\beta})^2
\le \pi ((2m)^\alpha + r (2m)^{\alpha\beta})^2
    - \pi (m^\alpha - r m^{\alpha\beta})^2
\]
for any $m \ge m_0$.
So $m \pi (r m^{\alpha\beta})^2
\le \pi ((2m)^\alpha + r (2m)^{\alpha\beta})^2$.
By taking $m \to \infty$, we obtain $1 + 2\alpha\beta \le 2\alpha$.

Next suppose that $-\infty <  \beta \le 0$.
For any $m > m_0$, we have
\[ \sum_{m_0 < n \le m} \pi (r n^{\alpha\beta})^2
\le \pi (m^\alpha + r m^{\alpha\beta})^2 - \pi (m_0^\alpha - r m_0^{\alpha\beta})^2,
\]
so
\[ 
(m - m_0) \pi (r m^{\alpha\beta})^2
\le  \pi (m^\alpha + r m^{\alpha\beta})^2.
\]
By taking $m\to\infty$, we have $1 + 2\alpha\beta \le 2\alpha$.
\end{proof}

\begin{proof}[Proof of Lemma \ref{43}]
Suppose that $\theta = p/q$ is an irreducible fraction.
Then for any $r>0$ and any $t > t_0 := r^{1/(1-\beta)}$,
we have $r t^\beta < t$ and
\begin{align*}
  B\left(t \rme^{\pi\rmi/q}, r t^\beta \sin\frac{\pi}{q} \right)
  \cap \Gamma_{\alpha,\theta} 
&\subset B\left(t \rme^{\pi\rmi/q}, t \sin\frac{\pi}{q} \right)
  \cap 
  \set{s \rme^{2\pi k\rmi / q} \mid s\ge0, k\in\zahl}
  \\
& = \emptyset.
\end{align*}
This implies that $\Gamma_{\alpha,\theta}$ is not asymptotically
$\beta$-relatively dense.

Let $z_n = n^\alpha \rme^{2\pi n\theta\rmi}$.
We have 
\[ \frac{|z_{(j+1)q}-z_{jq}|}{|z_{jq}|^\beta}
  = \frac{(j+1)^\alpha q^\alpha - (jq)^\alpha}{(jq)^{\alpha-\half}}
  = \sqrt{jq}\left( \left(1 + \frac{1}{j} \right)^\alpha - 1 \right)
  \to 0
\]
as $j \to \infty$.  This implies that $\Gamma_{\alpha,\theta}$ is not
asymptotically $\beta$-uniformly discrete.
\end{proof}

\section{Continued fractions and rational approximations}

This section prepares some properties of continued
fractions and rational approximations.
  A fraction $\frac{p}{q} \in \rat$ always assumes that
$p\in\zahl$ and $q \in \zahl_{>0}$. 
A pair of fractions $\frac{a}{m} < \frac{b}{n}$
is called a \textit{Farey pair} if $mb-na=1$.
 An open interval $(\frac{a}{m}, \frac{b}{n})$ 
 is called a \textit{Farey interval} of
if its endpoints $\frac{a}{m}, \frac{b}{n}$
are a Farey pair.

Let
\begin{equation}
\label{25}
 x
  = a_0+\frac{1}{a_1+\frac{1}{a_2+\cdots}} = [a_0,a_1,a_2,\cdots],
  \ a_0 \in \Z, \ a_i \in \Z_{>0}, \ i \in \zahl_{>0}
\end{equation}
be a continued fraction expansion of $x \in \real$.  
An irrational $x$ is called {\em badly approximable} if the set
of the partial quotients $\set{a_i \mid i \in \zahl_{\ge0}}$ is bounded.
Define the sequences $\set{p_i}_{i\geq -1}$
and $\set{q_i}_{i\geq -1}$ by $p_{-1} = 1$, $q_{-1} = 0$, $p_{0}=
a_0$, $q_{0}=1$, $p_{1}=a_0 a_1 + 1$, $q_{1}=a_1$, and $p_{i+1} =
a_{i+1} p_{i} + p_{i-1}$, $q_{i+1} = a_{i+1} q_{i} + q_{i-1}$,
$i\ge1$.  Let $p_{i,k} = k p_{i} + p_{i-1}$, $q_{i,k} = k q_{i} +
q_{i-1}$ for $i \ge0$, $0 \le k \le a_{i+1}$.  Note that $p_{i,0} =
p_{i-1}$, $q_{i,0} = q_{i-1}$, $p_{i,a_{i+1}} = p_{i+1}$,
$q_{i,a_{i+1}} = q_{i+1}$.  The fraction ${p_i}/{q_i} =
[a_0,a_1,\cdots,a_i]$, $i\ge0$, is called a {\it (principal)
  convergent} of $x$, and ${p_{i,k}}/{q_{i,k}} = [a_0, a_1, \cdots,
  a_i, k]$, $i\ge0$, $0<k<a_{i+1}$, is called an {\it intermediate
  convergent} of $x$.  An induction shows that if $i$ is odd and $0
\le k \le a_{i+1}$, then $(\frac{p_{i,k}}{q_{i,k}}, \frac{p_i}{q_i})$
is a Farey interval and $\frac{p_{i,k}}{q_{i,k}} < x <
\frac{p_i}{q_i}$.  If $i$ is even, $0 \le k \le a_{i+1}$ and
$(i,k)\neq(0,0)$, then $(\frac{p_i}{q_i}, \frac{p_{i,k}}{q_{i,k}})$ is
a Farey interval and $\frac{p_i}{q_i} < x < \frac{p_{i,k}}{q_{i,k}}$.

\begin{lemma} 
\label{38}
Let $(\frac{a}{m}, \frac{b}{n})$ be a Farey interval containing $x$.
Then we have either $(\frac{a}{m}, \frac{b}{n}) =(\frac{p_i}{q_i},
\frac{p_{i,k}}{q_{i,k}})$ for some $i$ even and $0 \le k < a_{i+1}$,
or $(\frac{a}{m}, \frac{b}{n}) =(\frac{p_{i,k}}{q_{i,k}},
\frac{p_i}{q_i})$ for some $i$ odd and $0 \le k < a_{i+1}$.
\end{lemma}
\begin{proof}
See \cite[Lemma 4]{arxiv2020}.
\end{proof}

\begin{lemma}
\label{8}
Let $\theta \in \real \setminus \rat$.
In the continued fraction expansion of $\theta$, the following inequalities hold.
\begin{enumerate}
\item
$q_{i,k} |q_{i,k}\theta - p_{i,k}| > \frac{1}{2+a_i}$
for $i\ge0$, $0 \le k < a_{i+1}$.
\item
$q_{i,k} |q_{i,k}\theta - p_{i,k}| < k+1$
for $i\ge0$, $0 \le k \le a_{i+1}$.
\item
$q_i |q_i \theta - p_i| < \frac{1}{a_{i+1}}$
for $i\ge0$.
\item
If $i\ge1$ and $k= \floor{\frac{1}{2}a_{i+1}}$,
then $q_{i,k} |q_{i,k}\theta - p_{i,k}| > \frac{a_{i+1} - 1}{4}$.
\end{enumerate}
\end{lemma}
\begin{proof}
1.
$|\theta - \frac{p_{i,k}}{q_{i,k}}|
> |\frac{p_{i+1}}{q_{i+1}} - \frac{p_{i,k}}{q_{i,k}}|
= \frac{a_{i+1}-k}{q_{i+1}q_{i,k}}$.
So $q_{i,k}|q_{i,k}\theta - p_{i,k}|
> \frac{q_{i,k}(a_{i+1}-k)}{q_{i+1}}
= \frac{(kq_i + q_{i-1})(a_{i+1}-k)}{a_{i+1}q_i+q_{i-1}}
\ge \frac{q_{i-1}a_{i+1}}{a_{i+1}q_i + q_{i-1}}
\ge \frac{q_{i-1}}{q_i+q_{i-1}}
> \frac{1}{2 + a_i}$.

2.
$|\theta - \frac{p_{i,k}}{q_{i,k}}|
< |\frac{p_{i}}{q_{i}} - \frac{p_{i,k}}{q_{i,k}}|
= \frac{1}{q_iq_{i,k}}$.
So $q_{i,k}|q_{i,k}\theta - p_{i,k}|
< \frac{q_{i,k}}{q_i} = \frac{kq_i + q_{i-1}}{q_i} \le k+1$.

3.
$|\theta - \frac{p_i}{q_i}| < |\frac{p_{i+1}}{q_{i+1}} - \frac{p_i}{q_i}|
= \frac{1}{q_i q_{i+1}}$.
So $q_i |q_i \theta - p_i| < \frac{q_i}{q_{i+1}} < \frac{1}{a_{i+1}}$.

4.
$|\theta - \frac{p_{i,k}}{q_{i,k}}|
> |\frac{p_{i+1}}{q_{i+1}} - \frac{p_{i,k}}{q_{i,k}}|
= \frac{a_{i+1}-k}{q_{i+1}q_{i,k}}$.
So $q_{i,k} |q_{i,k}\theta - p_{i,k}|
> \frac{q_{i,k}(a_{i+1}-k)}{q_{i+1}}
> \frac{k(a_{i+1}-k)}{a_{i+1}}$.
If $a_{i+1}$ is even, we have $k = \frac{a_{i+1}}{2}$, and
$\frac{k(a_{i+1}-k)}{a_{i+1}} = \frac{a_{i+1}}{4}$.
If $a_{i+1}$ is odd, we have $k = \frac{a_{i+1}-1}{2}$, and
$\frac{k(a_{i+1}-k)}{a_{i+1}} > \frac{a_{i+1}-1}{4}$.
\end{proof}

\section{Delone families of linear lattices}

This section considers relative denseness and uniform discreteness of a family of
(linear) lattices, and proves
Proposition \ref{7}.

Let $T \subset \real$ be a parameter set.
A family of point sets $\set{ \Lambda_t \subset \real^2 \mid t \in T}$ is called
relatively dense
if there exists $r>0$ such that for any $t\in T$ and any $\zeta \in \real^2$, we have $B(\zeta,r) \cap \Lambda_t \neq \emptyset$.
The family $\set{\Lambda_t}_{t \in T}$ is called uniformly discrete
if there exists $s>0$ such that for any $t \in T$ and any $\zeta \in \real^2$ we have $\#(B(\zeta,s) \cap \Lambda_t) \le 1$.
The family $\set{\Lambda_t}_{t}$ is called a Delone family
if it is relatively dense and uniformly discrete.

Let $\theta \in \real$, $z = (\theta, 1) \in \real^2$.
Consider the continued fraction expansion of $\theta$, as in the previous section.
Let $z_i = (q_i\theta - p_i, q_i)$ for $i\ge-1$,
and $z_{i,k} = (q_{i,k}\theta - p_{i,k}, q_{i,k})$ for $i\ge0$, $0 \le k \le a_{i+1}$.
We have $z_{-1} = z_{0,0} = (-1,0)$, $z_0 = (\theta - a_0, 1)$, $z_{0,1} = (\theta- a_0 - 1, 1)$, $z_{0,a_1} = z_1 = (q_1 \theta - p_1, q_1)$.
Let $\Lambda = z_{-1} \zahl + z_{0} \zahl
= \set{ (m\theta-n, m) \mid m,n\in\zahl}$ a lattice.
Let $T = \set{t \in \real \mid t \ge 1}$ from now on.
Let $A(t) = \begin{pmatrix} \sqrt{t} & 0 \\ 0 & 1/\sqrt{t} \end{pmatrix}$ be a $2\times2$ matrix, and
\[
\Lambda(t) := A(\sqrt{t}) \Lambda
 = \Set{ \left( (m\theta - n)\sqrt{t}, \frac{m}{\sqrt{t}} \right) \mid m,n\in\zahl}.
\]

Let
\begin{align*}
t_i &= \frac{q_i}{|q_i\theta - p_i|},\quad  i \ge -1, \\
t_{i,k} &= \frac{q_{i,k}}{|q_{i,k}\theta - p_{i,k}|},
\quad i\ge0, \ 0\le k \le a_{i+1}.
\end{align*}
We have $t_{-1} = t_{0,0} = 0$, $t_0 = \frac{1}{|\theta - a_0|}$,
$t_{0,1} = \frac{1}{|\theta - a_0 - 1|}$, $t_{0,a_1} = t_1 = \frac{1}{|\theta - a_0 - 1/a_1|}$,
and
\begin{gather*}
t_{i-1} < t_{i}, \\
t_{i,0} = t_{i-1} < t_{i,1} < \dots < t_{i,a_{i+1}} = t_{i+1}, \\
\sqrt{t_{i-1}t_i} = \sqrt{t_i t_{i,0}} 
< \sqrt{t_i t_{i,1}} < \dots < \sqrt{t_i t_{i,a_{i+1}}} = \sqrt{t_i t_{i+1}}
\end{gather*}
for $i \ge 0$.

Denote by
\begin{align*}
z_i(t)     &:= A(t) z_i = ((q_i\theta - p_i)\sqrt{t}, q_i/\sqrt{t}), \\
z_{i,k}(t) &:= A(t) z_{i,k} = ((q_{i,k}\theta - p_{i,k})\sqrt{t}, q_{i,k}/\sqrt{t}).
\end{align*}
We have 
\begin{align*}
z_i(t_i) &= ((-1)^i \sqrt{q_i|q_i\theta - p_i|}, \sqrt{q_i|q_i\theta - p_i|}), \\
z_{i,k}(t_{i,k}) &= ((-1)^{i+1} \sqrt{q_{i,k}|q_{i,k}\theta - p_{i,k}|}, \sqrt{q_{i,k}|q_{i,k}\theta - p_{i,k}|}).
\end{align*}

\begin{lemma}[Richards' formula]
Let $i\ge0$, $0\le k \le a_{i+1}$.  Let 
\[ t := \sqrt{t_i t_{i,k}}
  = \sqrt{\frac{q_i}{|q_i\theta - p_i|} 
\frac{q_{i,k}}{|q_{i,k}\theta - p_{i,k}|}}.
\] 
Then the parallelogram
$Q = \square(0, z_i(t), z_i(t) + z_{i,k}(t), z_{i,k}(t))$ is a rectangle. 
\end{lemma}
\begin{proof}
$z_i(t) = ((q_i\theta - p_i)\sqrt{t}, q_i/\sqrt{t})$,
$z_{i,k}(t) = ((q_{i,k}\theta - p_{i,k})\sqrt{t}, q_{i,k}/\sqrt{t})$,
and we have
\[ 
 (q_i\theta - p_i)\sqrt{t} \cdot (q_{i,k}\theta - p_{i,k})\sqrt{t}
 + \frac{q_i}{\sqrt{t}} \cdot \frac{q_{i,k}}{\sqrt{t}}
 = 0
\]
when $t = \sqrt{t_i t_{i,k}}$.
\end{proof}

\begin{lemma}
\label{22}
If $\set{\Lambda(t)}_{t\ge1}$ is relatively dense, then $\theta$ is badly approximable.
\end{lemma}
\begin{proof}
Suppose that $\sup_i a_i = \infty$.
Let $i \ge 1$, $k = \floor{\frac{a_{i+1}}{2}}$,
$r_i := \sqrt{q_{i,k} |q_{i,k}\theta - p_{i,k}|}$.
If $i$ is odd, let $Q(i) = [0, r_i] \times [0, r_i]$ be a square.
If $i$ is even, let $Q(i) = [- r_i, 0] \times [0, r_i]$. 
In either case, we have $\Lambda(t_{i,k}) \cap Q(i) = \set{0, z_{i,k}(t_{i,k})}$.
If $i$ is odd, we have $Q(i) \supset B(r_i(1+\rmi)/2, r_i/2)$
and $B(r_i(1+\rmi)/2, r_i/2) \cap \Lambda(t_{i,k}) = \emptyset$.
If $i$ is even, we have $Q(i) \supset B(r_i(-1+\rmi)/2, r_i/2)$
and  $B(r_i(-1+\rmi)/2, r_i/2) \cap \Lambda(t_{i,k}) = \emptyset$.
By Lemma~\ref{8}, we have $\sup_i r_i \ge \sup_i \sqrt{a_{i+1}-1}/2 = \infty$,
so $\set{\Lambda(t)}_{t\ge1}$ is not relatively dense.
\end{proof}

\begin{lemma}
\label{23}
If $\set{\Lambda(t)}_{t\ge1}$ is uniformly discrete, then $\theta$ is badly approximable.
\end{lemma}
\begin{proof}
Suppose that $\sup_i a_i = \infty$.
We have $\dist(z_i(t_i),0) = |z_i(t_i)| =  \sqrt{2 q_{i} |q_{i}\theta - p_{i}|} < \sqrt{2/a_{i+1}}$,
so $\inf_i \dist(z_{i}(t_i),0) = 0$,
and the family $\set{\Lambda(t)}_{t\ge1}$ is not uniformly discrete.
\end{proof}

\begin{lemma}
\label{24}
If $\theta$ is badly approximable, then $\set{\Lambda(t)}_{t\ge1}$ is uniformly discrete.
\end{lemma}
\begin{proof}
Suppose that $\sup_i a_i = M < \infty$.
We have 
\begin{align*}
|z_{i,k}(t)|^2
&= t(q_{i,k}\theta - p_{i,k})^2 + q_{i,k}^2/t \\
&\ge 2 q_{i,k}|q_{i,k}\theta - p_{i,k}|
\ge \frac{2}{2+a_i} \ge \frac{2}{2+M}.
\end{align*}
This implies that $\dist(\lambda, 0) = |\lambda| \ge \sqrt{\frac{2}{2+M}}$ for any 
 $t \ge1$, $\lambda \in \Lambda(t) \setminus\set{0}$,
and that $\dist(\lambda, \lambda') = |\lambda - \lambda'| \ge \sqrt{\frac{2}{2+M}}$ for any distinct $\lambda, \lambda' \in \Lambda(t)$.
So $\set{\Lambda(t)}_{t\ge1}$ is uniformly discrete.
\end{proof}

\begin{lemma}
\label{5}
If $\theta$ is badly approximable,
then $\set{\Lambda(t)}_{t\ge1}$ is relatively dense.
\end{lemma}
\begin{proof}
Let $M = \sup_i a_i < \infty$.
Fix $t \ge 1$.
There exist $i\ge -1$ such that
$t_i \le t \le t_{i+1}$.
Let
\[ 
  Q = \square(0,z_i(t), z_i(t) + z_{i+1}(t), z_{i+1}(t))
\]
be a closed parallelogram.
The family of translated parallelograms $\lambda +Q$, $\lambda \in \Lambda(t)$, covers the plane,
$\bigcup_{\lambda \in \Lambda(t)} (\lambda + Q) = \cx$.
For any $\zeta \in \real^2$, there exists $\lambda \in \Lambda(t)$ such that
$\zeta \in \lambda + Q$.
Let
\[
 Q' = [- |\fracpart{q_i\theta}|\sqrt{t}, |\fracpart{q_i\theta}|\sqrt{t}]
  \times \left[0, \frac{q_i + q_{i+1}}{\sqrt{t}} \right] \supset Q
\]
be a rectangle.  Since $\zeta \in \lambda + Q'$,
we have
\[ |\zeta - \lambda| \le |\fracpart{q_i\theta}|\sqrt{t} + \frac{q_i + q_{i+1}}{\sqrt{t}},
\]
where
\[
|q_{i}\theta - p_{i}| \sqrt{t}
\le |q_{i}\theta - p_{i}| \sqrt{t_{i+1}} 
= \sqrt{\frac{q_{i+1}(q_i\theta - p_i)^2}{|q_{i+1}\theta - p_{i+1}|}} 
\le \sqrt{a_{i+2}+1} \le \sqrt{M+1},
\]
\begin{align*}
\frac{q_i + q_{i+1}}{\sqrt{t}} 
&\le \frac{q_i + q_{i+1}}{\sqrt{t_i}} 
= \sqrt{\frac{(q_i+q_{i+1})^2 |q_i\theta - p_i|}{q_i}} \\
&\le \sqrt{\frac{(q_i+q_{i+1})^2}{q_i q_{i+1}}} 
\le \sqrt{\frac{(M+2)^2}{M+1}}
\le 2 \sqrt{M+1}.
\end{align*}
Thus, for any $\zeta \in \real^2$, there exist $\lambda \in \Lambda(t)$ such that
$|\zeta - \lambda| \le 3 \sqrt{M+1}$.
So $\set{\Lambda(t)}_{t\ge1}$ is relatively dense.
\end{proof}

\begin{proof}[Proof of Proposition~\ref{7}]
The proof is given by Lemmas \ref{22}-\ref{5}.
\end{proof}

\section{Spiral Delone set}

This section proves Theorem \ref{44}.
Let $\alpha >0$, $- \infty < \beta = 1 - \frac{1}{2\alpha} < 1$.
Let
$\Gamma_{\alpha,\theta} = \set{F(n) \mid n \in \zahl_{>0}}$,
where $F(n) := n^\alpha \rme^{2\pi n \theta \rmi}$.

\begin{lemma}
\label{12}
Let $c>0$, $0<\alpha<1$.
If $0<x\le 1$, then we have
\[
 (1+cx)^\alpha - 1 \ge x((1+c)^\alpha - 1).
\]
\end{lemma}
\begin{proof}
Let $f(x) = (1+cx)^\alpha$.  Since the function $-f(x)$ is convex, 
we have 
\[
\frac{f(x) - f(0)}{x-0} = \frac{(1+cx)^\alpha - 1}{x}  
\ge \frac{f(1)-f(0)}{1-0}
= (1+c)^\alpha - 1
\]
for $0<x \le 1$.
\end{proof}

\begin{lemma}
\label{21}
If $0 < \alpha < 1$
and $0 < x \le \frac{1}{1-\alpha}$, 
then $(1+x)^\alpha \ge 1 + \frac{\alpha x}{2}$.
\end{lemma}
\begin{proof}
By Taylor's theorem, there exists $0 < \xi < x$ such that 
$(1+x)^\alpha = 1 + \alpha x + \frac{\alpha(\alpha-1)}{2}(1+\xi)^{\alpha-2}x^2$.
So we have
\begin{align*}
(1+x)^\alpha &= 1 + \alpha x + \frac{\alpha(\alpha-1)}{2}(1+\xi)^{\alpha-2}x^2 \\
&> 1 + \alpha x + \frac{\alpha(\alpha-1)}{2} x^2 \\
&\ge 1 + \frac{\alpha x}{2} 
\end{align*}
for $x \le 1/(1-\alpha)$.
\end{proof}

\begin{lemma}
\label{45}
Suppose that $\theta$ is badly approximable.
Then $\Gamma_{\alpha,\theta}$ is asymptotically $\beta$-uniformly discrete.
\end{lemma}
\begin{proof}
Suppose that $\sup_i a_i = M < \infty$.
By Lemma~\ref{8}, we have $q_{i,k}|q_{i,k}\theta - p_{i,k}| > \frac{1}{2+a_i} \ge \frac{1}{2+M}$ for $i\ge0$, $0\le k < a_{i+1}$.
This implies that
 $\frac{q}{\sqrt{\nu}}\cdot |q\theta -p|\sqrt{\nu} > \frac{1}{2+M}$
for any $\nu, q \in \zahl_{>0}$ and $p\in\zahl$.

Fix $\nu, q \in \zahl_{>0}$.  Then we have either
\begin{equation}
\label{11}
|\fracpart{q\theta}|\sqrt{\nu} \ge \frac{1}{\sqrt{2+M}}
\end{equation}
or
\begin{equation}
\label{10}
\frac{q}{\sqrt{\nu}} \ge \frac{1}{\sqrt{2+M}}.
\end{equation}
Suppose first that (\ref{11}) holds.
We have
\begin{align*}
|F(\nu+q) - F(\nu)| 
&\ge \nu^\alpha |\rme^{2\pi \theta q \rmi} - 1| \\
&\ge 2\nu^\alpha \sin|\fracpart{q\theta}\pi| 
\ge 4\nu^\alpha |\fracpart{q\theta}| \\
&\ge |F(\nu)|^{\beta} \frac{4}{\sqrt{2+M}}.
\end{align*}
Next suppose that (\ref{10}) holds.
If $\alpha \ge 1$,
we have
\begin{align*}
|F(\nu+q) - F(\nu)| 
&\ge (\nu + q)^\alpha - \nu^\alpha
= |F(\nu)|^\beta \sqrt{\nu} \left((1+\frac{q}{\nu})^\alpha - 1 \right) \\
&\ge |F(\nu)|^\beta \sqrt{\nu} \left( (1+\frac{q}{\nu}) - 1 \right)
= |F(\nu)|^\beta \cdot \frac{q}{\sqrt{\nu}} \\
&\ge |F(\nu)|^\beta \frac{1}{\sqrt{2+M}}.
\end{align*}
If $0 < \alpha < 1$, 
we have
\begin{align*}
|F(\nu+q) - F(\nu)| 
&\ge (\nu + q)^\alpha - \nu^\alpha
= \nu^\alpha \left((1+\frac{q}{\nu})^\alpha - 1 \right) \\
&\ge \nu^{\alpha} 
 \left( \left( 1+\frac{1}{\sqrt{\nu(2+M)}} \right)^\alpha - 1 \right) \\
&\ge 
\nu^{\alpha - \half} 
 \frac{2^\alpha - 1}{\sqrt{2+M}}
= |F(\nu)|^{\beta} 
 \frac{2^\alpha - 1}{\sqrt{2+M}}
\end{align*}
by Lemma~\ref{12}.
Thus we obtain
\[
|F(\nu +q) - F(\nu)|
\ge |F(\nu)|^\beta  \frac{\min\set{1,2^\alpha - 1}}{\sqrt{2+M}} .
\]
for any $\nu,q\in\zahl_{>0}$.
This implies that $\Gamma_{\alpha,\theta}$ is $\beta$-uniformly discrete.
\end{proof}

\begin{lemma}
If $\Gamma_{\alpha,\theta}$ is asymptotically $\beta$-uniformly discrete, then $\theta$ is badly approximable.
\end{lemma}
\begin{proof}
Let $\nu_i := \ceil{t_i} \in \zahl$. 
We have
\begin{align*}
|F(\nu_i + q_i) - F(\nu_i)|
&\le |F(\nu_i + q_i) - \nu_i^\alpha \rme^{2\pi \theta (\nu_i + q_i) \rmi}|
 + |\nu_i^\alpha \rme^{2\pi \theta (\nu_i + q_i) \rmi} - F(\nu_i)| \\
&\le ((\nu_i + q_i)^\alpha - \nu_i^\alpha) + \nu_i^\alpha 2\pi |\fracpart{q\theta}| \\
&= |F(\nu_i)|^\beta \left(
 \sqrt{\nu_i} ((1+\frac{q_i}{\nu_i})^\alpha - 1) + 2\pi |\fracpart{q\theta}| \sqrt{\nu_i}
\right)
\end{align*}
where
\[
|\fracpart{q_i\theta}| \sqrt{\nu_i}
< |\fracpart{q_i\theta}| \sqrt{t_i + 1} 
<  |\fracpart{q_i\theta}| \sqrt{2 t_i} 
= \sqrt{2 q_i|q_i\theta - p_i|} 
< \sqrt{\frac{2}{a_{i+1}}}.
\]
If $i$ is large, we may assume that $\frac{q_i}{\nu_i} \le \frac{q_i}{t_i} = |q_i\theta - {p_i}| \le \frac{1}{\alpha^2}$, which implies
\begin{align*}
\sqrt{\nu_i} ((1+\frac{q_i}{\nu_i})^\alpha - 1)
&\le \sqrt{\nu_i} (1+\frac{2q_i}{\nu_i} - 1)
 \le \frac{2 q_i}{\sqrt{t_i}}
 = 2 \sqrt{q_i|q_i\theta - p_i|}
  < \frac{2}{\sqrt{a_{i+1}}} .
\end{align*}
Thus we obtain
\begin{align*}
|F(\nu_i + q_i) - F(\nu_i)|
\le \frac{2+\sqrt{2}}{\sqrt{a_{i+1}}} |F(\nu_i)|^\beta 
\end{align*}
for $i$ large.
Since $\Gamma_{\alpha,\theta}$ is $\beta$-uniformly discrete, we have $\sup_i a_i < +\infty$.
\end{proof}

\begin{lemma}
\label{13}
Suppose that $\theta$ is badly approximable.
Then $\Gamma_{\alpha,\theta}$ is asymptotically $\beta$-relatively dense.
\end{lemma}
\begin{proof}
Consider the polar coordinates of the plane, $\varphi : \real_{>0} \times \real \to \cx$,
$\varphi(s, t) = s^\alpha \rme^{2\pi t \rmi}$.
We have 
$\varphi(k,k\theta) = F(k)$
and
$\varphi(s,t+k) = \varphi(s,t)$  
 for $k \in \zahl$.
Let $\Lambda = (1,\theta)\zahl + (0,1)\zahl$ be a linear lattice.
We have $\varphi(\Lambda) = \Gamma_{\alpha,\theta}$.

Suppose $i_0$ is sufficiently large that 
\begin{equation}
\label{14}
q_{i_0} \ge 2(1+\alpha^2), 
\qquad 
q_{i_0} \ge 3(2+M).
\end{equation}
Let $\zeta \in \cx$.  
We are going to show that
if $|\zeta| \ge t_{i_0}$, then
there exists $\nu \in \zahl_{>0}$ such that
\begin{equation}
\label{20}
|\zeta - F(\nu)| \le |F(\nu)|^\beta ( 2\alpha\sqrt{2(M+2)} + 2\pi\sqrt{M+1}).
\end{equation}
If $|\zeta| \ge {t_{i_0}}$, there exists $i \ge i_0$ such that
${t_i} \le |\zeta| < {t_{i+1}}$.
Let 
\[ T = \square((0,0), (q_i, \fracpart{q_i\theta}), (q_i+q_{i+1}, \fracpart{(q_i+q_{i+1})\theta}), (q_{i+1}, \fracpart{q_{i+1}\theta}))
\]
be a parallelogram.
There exists $\nu \in \zahl_{>0}$ such that
$\zeta \in \varphi((\nu, \fracpart{\nu\theta}) + T)$.
We have ${t_i} - (q_i + q_{i+1}) \le \nu \le {t_{i+1}}$.
There exist $s,t\in\real$ such that
$\zeta = (\nu+s)^\alpha \rme^{2\pi(\fracpart{\nu\theta} + t)\rmi}$,
$0\le s \le q_i + q_{i+1}$, $|t| \le |\fracpart{q_i\theta}|$.
We have
\begin{align*}
& |\zeta - F(\nu)| \\
&= |(\nu+s)^\alpha \rme^{2\pi (\nu\theta+t) \rmi} - \nu^\alpha \rme^{2\pi\nu\theta \rmi}| \\
&\le 
|(\nu+s)^\alpha \rme^{2\pi(\nu\theta+t) \rmi} - \nu^\alpha \rme^{2\pi(\nu\theta+t)\rmi}| 
+ |\nu^\alpha \rme^{2\pi (\nu\theta+t) \rmi} - \nu^\alpha \rme^{2\pi\nu\theta\rmi}| \\
&\le ((\nu+s)^\alpha - \nu^\alpha)
+ \nu^\alpha 2\pi |t| \\
&= |F(\nu)|^\beta \left( \sqrt{\nu} \left( \left(1+\frac{s}{\nu} \right)^\alpha - 1  \right)
 + 2\pi |t| \sqrt{\nu} \right) .
\end{align*}
We have
\begin{align*}
2\pi |t| \sqrt{\nu}
&\le 2\pi |q_i\theta - p_i| \sqrt{t_{i+1}} 
= 2\pi \sqrt{ \frac{q_{i+1} (q_i\theta - p_i)^2}{|q_{i+1}\theta - p_{i+1}|} } \\
&\le 2\pi \sqrt{a_{i+2}+1} \\
&\le 2\pi \sqrt{M+1}.
\end{align*}
By (\ref{14}), we have
\begin{align*}
\frac{\nu}{s} + 1 
&\ge \frac{{t_i}-(q_i+q_{i+1})}{q_i + q_{i+1}} + 1 
= \frac{t_i}{q_i + q_{i+1}}\\
&=  {\frac{q_i}{(q_i + q_{i+1}) |q_i\theta - p_i|}} \\
&\ge  \frac{q_i q_{i+1}}{q_i + q_{i+1}} 
 \ge \frac{q_i}{2} 
 \ge 1+\alpha^2,
\end{align*}
so $0 < \frac{s}{\nu} \le \frac{1}{\alpha^2}$.
By Lemma~\ref{18}, we have
\[
\sqrt{\nu} \left( \left(1+\frac{s}{\nu} \right)^\alpha - 1  \right)
\le \sqrt{\nu} \left( 1 + \frac{2\alpha s}{\nu} - 1 \right) \\
=  \frac{2\alpha s}{\sqrt{\nu}}.
\]
We have
\begin{align*}
\frac{\sqrt{\nu}}{s}
&\ge \frac{\sqrt{t_i - (q_i+q_{i+1})}}{q_i + q_{i+1}} 
 = \sqrt{\frac{q_iq_{i+1}}{(q_i + q_{i+1})^2} - \frac{1}{q_i+q_{i+1}}} \\
&\ge \sqrt{\frac{M+1}{(M+2)^2} - \frac{1}{2q_i}} \\
&\ge  \sqrt{\frac{M+1}{(M+2)^2} - \frac{1}{6(M+2)}} 
= \sqrt{\frac{5M+4}{6(M+2)^2}} \\
&= \frac{1}{\sqrt{2(M+2)}} .
\end{align*}
Thus we obtain (\ref{20}).
\end{proof}

\begin{lemma}
\label{46}
If $\Gamma_{\alpha,\theta}$ is asymptotically $\beta$-relatively dense, then $\theta$ is badly approximable.
\end{lemma}
\begin{proof}
Suppose that $\sup_i a_i = \infty$.
Let $i\ge0$, $k = \floor{\frac{a_{i+1}}{2}}$,
$\nu_i := \ceil{t_{i,k}}$.
Let
\[
 D := \set{ x \rme^{2\pi y\rmi} \mid \nu_i^\alpha \le x \le (\nu_i + q_{i,k})^\alpha,
0 \le y - \nu_i\theta \le q_{i,k}\theta - p_{i,k}}
\]
if $i$ is odd, or
\[
 D := \set{ x \rme^{2\pi y\rmi} \mid \nu_i^\alpha \le x \le (\nu_i + q_{i,k})^\alpha,
0 \ge y - \nu_i\theta \ge q_{i,k}\theta - p_{i,k}}
\]
if $i$ is even.
In either case we have 
$D \cap \Gamma_{\alpha,\theta} 
= \set{ F(\nu_i), F(\nu_i+q_{i,k})}$,
and
$D \supset B(\xi_i, R_i)$, 
where
\[ \xi_i := \frac{\nu_i^\alpha + (\nu_i + q_{i,k})^\alpha}{2} 
          \rme^{2\pi(\nu_i\theta + (q_{i,k}\theta - p_{i,k})/2) \rmi},
\]
\[
 R_i := \min\Set{\frac{(\nu_i + q_{i,k})^\alpha - \nu^\alpha}{2}, 
   |\xi_i| \sin|\pi(q_{i,k}\theta - p_{i,k})|}.
\]
So we have
\[
  |\xi_i - F(m)| \ge |\xi_i|^\beta r_i
\]
for any $m \in \zahl_{>0}$, where $r_i := R_i |\xi_i|^{-\beta}$.
If $\Gamma_{\alpha,\theta}$ is $\beta$-relatively dense,
then $\sup_i r_i < +\infty$.

Since
$|\xi_i|^{1/\alpha} > \nu \ge t_{i,k}$, we have
\begin{align*}
|\xi_i| \sin|\pi(q_{i,k}\theta - p_{i,k})| 
&\ge 2 |\xi_i| |q_{i,k}\theta - p_{i,k}| \\
&\ge 2 |\xi_i|^{\beta} \sqrt{t_{i,k}} |q_{i,k}\theta - p_{i,k}| \\
&= 2 |\xi_i|^{\beta} \sqrt{q_{i,k} |q_{i,k}\theta - p_{i,k}| } \\
&\ge |\xi_i|^{\beta} \sqrt{{a_{i+1}-1}}.
\end{align*}
We may assume that $i$ is so large that 
$\frac{q_{i,k}}{t_{i,k}} = |q_{i,k}\theta - p_{i,k}| \le 3^{1/\alpha} - 1$ holds, which imply that
\[ 
1 < \frac{|\xi_i|}{\nu^\alpha}
= \half \left(1 + \left( 1+\frac{q_{i,k}}{\nu} \right)^\alpha \right)
\le \half \left(1 + \left( 1+\frac{q_{i,k}}{t_{i,k}} \right)^\alpha \right)
 \le 2.
\]
If $\alpha \ge 1$, we have
\begin{align*}
(\nu + q_{i,k})^\alpha - \nu^\alpha
&= {\nu^\alpha} ((1 + \frac{q_{i,k}}{\nu})^\alpha - 1) \\
&\ge \frac{1}{2} |\xi_i| (1 + \frac{\alpha q_{i,k}}{\nu} - 1) \\
& > \frac{1}{2} |\xi_i|^{\beta} \frac{\alpha q_{i,k}}{\sqrt{\nu}} 
> \frac{1}{2} |\xi_i|^{\beta} \frac{\alpha q_{i,k}}{\sqrt{1 + \nu_{i,k}}} 
> \frac{1}{2} |\xi_i|^{\beta} \frac{\alpha q_{i,k}}{\sqrt{2\nu_{i,k}}} \\
&\ge \frac{\alpha}{2\sqrt{2}} |\xi_i|^{\beta} \sqrt{q_{i,k}|q_{i,k}\theta - p_{i,k}|} \\
&\ge \frac{\alpha}{4\sqrt{2}} |\xi_i|^{\beta} \sqrt{{a_{i+1}-1}}.
\end{align*}
If $0 < \alpha < 1$, we further assume that
 $\frac{q_{i,k}}{\nu} \le  \frac{q_{i,k}}{t_{i,k}} = |q_{i,k}\theta - p_{i,k}| \le \frac{1}{1-\alpha}$.
 By Lemma~\ref{21}, 
we have
\begin{align*}
(\nu + q_{i,k})^\alpha - \nu^\alpha
&= {\nu^\alpha} ((1 + \frac{q_{i,k}}{\nu})^\alpha - 1) \\
&\ge \frac{1}{2} |\xi_i|^{\beta} (1 + \frac{\alpha q_{i,k}}{2\nu} - 1) \\
&> \frac{\alpha}{4} |\xi_i|^{\beta} \frac{q_{i,k}}{\sqrt{\nu}} 
> \frac{\alpha}{4} |\xi_i|^{\beta} \frac{q_{i,k}}{\sqrt{2\nu_{i,k}}} \\
&\ge \frac{\alpha}{4\sqrt{2}} |\xi_i|^{\beta} \sqrt{q_{i,k}|q_{i,k}\theta - p_{i,k}|} \\
&\ge \frac{\alpha}{8\sqrt{2}} |\xi_i|^{\beta} \sqrt{{a_{i+1}-1}}.
\end{align*}
Thus we obtain $r_i \ge \min\set{1, \frac{\alpha}{16\sqrt{2}}} \sqrt{a_{i+1}-1}$.
If $\Gamma_{\alpha,\theta}$ is $\beta$-relatively dense,
then we have $\sup_i r_i < +\infty$ and $\sup_i a_i < +\infty$. 
\end{proof}

\begin{proof}[Proof of Theorem \ref{44}]
The proof is given by Lemmas \ref{45}-\ref{46}.
\end{proof}



\end{document}